\documentclass{article}
\usepackage[centering,total={410pt,610pt}]{geometry}
\usepackage{graphicx}
\usepackage[titletoc]{appendix}
\usepackage{tikz}
\usetikzlibrary{calc}
\usepackage{pgfplots}
\pgfplotsset{compat=newest}
\usepackage{appendix}
\usepackage{tikz}
\usepackage{enumerate}
\usepackage{amsmath}
\usepackage{amssymb}
\usepackage{amsthm}

\theoremstyle{definition}

\theoremstyle{remark}

\usepackage{hyperref}
\hypersetup{colorlinks,citecolor=black,filecolor=black,linkcolor=black,urlcolor=black}

\numberwithin{equation}{section}
\allowdisplaybreaks

\usepackage{color}

\newcommand{\R}{\mathbb{R}}

\title{On the dual step length of\break the alternating direction method of multipliers}

\author{%
Guoyong Gu%
\thanks{Department of Mathematics, Nanjing University, Nanjing, 210093, China.}
\thanks{Email: ggu@nju.edu.cn. This author was supported by the NSFC grant 11671195.}
\and
Junfeng Yang\footnotemark[1]
\thanks{Email: jfyang@nju.edu.cn.
This author was supported by the NSFC grants 11771208 and 11922111.}
}
\date{}

\begin{document}
\maketitle

\begin{abstract}
The alternating direction method of multipliers (ADMM) is a most widely used optimization scheme for solving linearly constrained separable convex optimization problems. The convergence of the ADMM can be guaranteed when the dual step length is less than the golden ratio, while plenty of numerical evidence suggests that even larger dual step length often accelerates the convergence. It has also been proved that the dual step length can be enlarged to less than 2 in some special cases, namely, one of the separable functions in the objective function is linear, or both are quadratic plus some additional assumptions. However, it remains unclear whether the golden ratio can be exceeded in the general convex setting.
In this paper, the performance estimation framework is used to analyze the convergence of the ADMM, and assisted by numerical and symbolic computations, a counter example is constructed, which indicates that the conventional measure may lose monotonicity as the dual step length exceeds the golden ratio, ruling out the possibility of breaking the golden ratio within this conventional analytic framework.

\bigskip

\noindent\textbf{Keywords:} Alternating direction method of multipliers, golden ratio, performance estimation framework, semidefinite programming.
\end{abstract}


\section{Introduction}
We consider linearly constrained separable convex optimization problem of the form
\begin{equation}\label{prob}
\min\bigl\{f(x)+g(y) \mid Ax+By=b, \; x\in\R^{n_1}, \; y\in\R^{n_2} \bigr\},
\end{equation}
where $f: \R^{n_1} \rightarrow \R \cup \{+\infty\}$ and $g: \R^{n_2} \rightarrow \R \cup \{+\infty\}$ are extended real-valued closed proper convex functions,
$A\in\R^{m\times n_1}$, $B\in\R^{m\times n_2}$ and $b\in\R^{m}$.
Throughout this paper, we let $\partial f$ and $\partial g$ be the subdifferential operators of $f$ and $g$, respectively,
and let $\|\cdot\|$ be the $\ell_2$ norm induced by the dot inner product $\langle u, v\rangle := u^Tv$, where $u$ and $v$ are column vectors of the same length.

We assume that the solution set of \eqref{prob} is nonempty.
Furthermore, constraint qualification holds so that $(x^*,y^*)\in\R^{n_1}\times\R^{n_2}$ is a solution of \eqref{prob} if and only if there exists  $z^*\in\R^{m}$ such that $(x^*,y^*,z^*)$ is a KKT point, i.e., $(x^*,y^*,z^*)\in\Omega^*$, where
\begin{equation}
  \label{kkt}
\Omega^* := \{(x^*,y^*,z^*) \in\R^{n_1}\times\R^{n_2}\times \R^m \mid A^Tz^* \in \partial f(x^*), \; B^Tz^* \in \partial g(y^*),\; Ax^*+By^*=b\}.
\end{equation}

Let the augmented Lagrangian function of \eqref{prob} be given by
\begin{equation*}
 L_\beta(x, y, z)=f(x)+g(y)-z^T(Ax+By-b)+\frac{\beta}{2}\|Ax+By-b\|^2,
\end{equation*}
where $\beta>0$ is the penalty parameter and $z\in\R^m$ is the dual variable (or Lagrange multiplier).
The alternating direction method of multipliers (ADMM) was originally proposed in \cite{GloM75,GabM76} for solving elliptic partial differential equations and nonlinear variational problems.
It solves \eqref{prob} via
\begin{equation}\label{admm}
\left\{
\begin{array}{l}
  x^{k+1}=\arg \min_{x}L_\beta(x, y^k, z^k),\\[0.2em]
  y^{k+1}=\arg \min_{y}L_\beta(x^{k+1}, y, z^k),\\[0.2em]
  z^{k+1}=z^k-\gamma\beta(Ax^{k+1}+By^{k+1}-b),
\end{array}
\right.
\end{equation}
for $k=0,1,2,\ldots$,  where $y^0\in\R^{n_1}$ and $z^0\in\R^m$ are given initial points and $\gamma >0$ is a dual step length.
Equivalently, the iterative scheme \eqref{admm} can be restated as
\begin{equation}\label{admm2}
\left\{
\begin{array}{l}
  A^T[z^k - \beta(Ax^{k+1}+By^k-b)] \in  \partial f(x^{k+1}),\\[0.2em]
  B^T[z^k - \beta(Ax^{k+1}+By^{k+1}-b)]  \in  \partial g(y^{k+1}),\\[0.2em]
  z^{k+1}=z^k-\gamma\beta(Ax^{k+1}+By^{k+1}-b).
\end{array}
\right.
\end{equation}
Due to its efficiency, scalability, versatility and robustness, the ADMM has become an important algorithmic tool for solving modern convex  optimization problems with desirable structures.
During the last decade, the ADMM has achieved enormous popularity in solving a large number of well structured
optimization problems arising from various applications including signal and image processing, engineering,  statistics, machine learning and big data analysis, and so on.
We refer to \cite{Boyd+11} and the references thereof for a sample of modern applications of the ADMM.

The convergence of the ADMM has been studied ever since it was invented \cite{GloM75,GabM76}.
Nowadays, the convergence and rate of convergence of ADMM are well understood from both perspectives of the augmented Lagrangian and the fixed point theory of non-expansive operators due to the connections to the Douglas-Rachford operator splitting method \cite{Gabay83} and the proximal point method \cite{EB92}.
Interested reader is referred to \cite{EY15} for a recent survey of ADMM from the perspective of fixed-point theory of nonexpansive operators and \cite{Glow14} for a historical perspective on ADMM.

The focus of this paper is about the dual step length $\gamma$, which plays an important role in practical computations.
It was firstly shown by Fortin and Glowinski \cite{FG83} that under some assumptions the ADMM converges for $\gamma \in (0, (\sqrt{5}+1)/2)$,
while numerical evidence has been observed in \cite{WGY10mpc} showing that even larger values of $\gamma$ usually lead to faster convergence.
Therefore, it is of both theoretical and practical importance to study the possibility of breaking the restriction of the golden ratio on $\gamma$.
Indeed, this has been pointed out as an open question in \cite[Remark 5.1]{Glow84}.
In the pioneering work \cite{GabM76}, Gabay and Mercier showed that the ADMM converges for the broader region $\gamma \in (0, 2)$ when $g$ is linear.
When both $f$ and $g$ are convex quadratic, it was recently shown in \cite{TY18jota} that ADMM converges for $\gamma \in (0, 2)$ if some additional conditions are satisfied.
However, it remains unclear if the golden ratio can be exceeded in the general convex setting.

Let $(x^*, y^*, z^*)\in\Omega^*$ be a KKT point of \eqref{prob}.
A conventional measure for establishing convergence of ADMM in the case $\gamma \in (1,(\sqrt{5}+1)/2))$ is
\begin{equation}\label{Rk}
  R^k := \|z^k - z^*\|^2 + \gamma \beta^2 \|B(y^k - y^*)\|^2 + (\gamma-1)\beta^2\|Ax^k+By^k-b\|^2.
\end{equation}
In fact, Fortin and Glowinski \cite[Eqs. (5.26) and (5.27)]{FG83} essentially showed that $R^k$ is monotonically decreasing when $(x^k,y^k,z^k)\notin \Omega^*$, which, together with some additional conditions, leads to the convergence of ADMM for $\gamma \in (1,(\sqrt{5}+1)/2))$. A similar analysis is given in \cite[Theorem 1]{HY98orl} and an extension to the semi-proximal ADMM can be found in \cite[Theorem B.1]{FPST13sima}.
In this paper, we study the possibility of breaking the golden ratio restriction on $\gamma$ using the idea of performance estimation proposed by Drori and Teboulle \cite{DT14mp,DT16mp}.
Assisted by numerical and symbolic computations, we show that in the worst case the measure $R^k$ defined in \eqref{Rk} fails to be monotonically decreasing for $\gamma > (\sqrt{5}+1)/2$, ruling out the possibility of extending $\gamma$ beyond the golden ratio within this analytic framework.

The paper is organized as follows. We formulate in section \ref{sc:pep} the performance estimation problem,
which is then equivalently reformulated in section \ref{sc:sdp} as a semidefinite programming problem (SDP).
A particular rank two feasible solution of the SDP is constructed in section \ref{sc:rank2},
which indicates that the conventional measure $R^k$ defined in \eqref{Rk} fails to be monotonically decreasing for $\gamma>(\sqrt{5}+1)/2$.
Finally, some concluding remarks are given in section \ref{sc:remarks}.

\section{Performance estimation}\label{sc:pep}
In the performance estimation framework,
the tight worst case performance of a specific algorithm designed for solving a particular class of problems is first represented as an infinite dimensional nonconvex optimization problem, which is then reduced to a finite dimensional nonconvex optimization problem by using ceratin interpolation conditions, and finally an equivalent SDP is formulated by a semidefinite lifting technique.
Pioneered by Drori and Teboulle \cite{DT14mp,DT16mp}, the performance estimation idea has recently attracted a significant amount of attention.
It has been applied to study the worst case convergence rates of various optimization algorithms, see the recent works by Kim and Fessler ~\cite{KF16MPA,KF18SIOPTa} and Taylor and his collaborators \cite{THG17mp,THG17siopt}.
Besides deterministic algorithms, the performance estimation idea has also been extended in \cite{TB19} to analyze stochastic first order algorithms.
Furthermore, it has been extended by Ryu et al. \cite{RTBG18} to study operator splitting methods for monotone inclusion problems.
See also our recent analysis of the proximal point algorithm within this framework \cite{GY19nonergodic}.

Let the measure $R^k$ be defined in \eqref{Rk} and $(x^{k+1},y^{k+1},z^{k+1})$ be generated via the ADMM scheme \eqref{admm} or \eqref{admm2} from $(x^k,y^k,z^k)$.
The worst case value of $R^{k+1}$ is represented as the following performance estimation problem (PEP)
\begin{equation}
  \label{pep}
\sup \left\{R^{k+1} \left|
  \begin{array}{c}
  f \text{~and~} g \text{~are closed proper convex}, (x^*,y^*,z^*)\in\Omega^*, \; R^k = 1,  \\[0.2em]
      (x^{k+1},y^{k+1},z^{k+1}) \text{~~is generated via \eqref{admm} from~~}(x^k,y^k,z^k)
  \end{array}
  \right.\right\}.
\end{equation}
Here, the optimization variables are $(f,g,A,B,b)$, $(x^*,y^*,z^*)$, $(x^k,y^k,z^k)$, $(x^{k+1},y^{k+1},z^{k+1})$ and $\beta$.
The constraint $R^k=1$ is a normalization condition.
It is clear that \eqref{pep} is infinite dimensional, nonconvex and too complicated to be solved directly.
The attractive feature of performance estimation is that \eqref{pep} can be equivalently reformulated as a finite dimensional convex SDP through a series of reductions and transformations, which we explain next.

\subsection{Invariance  with $\beta$}
Assume that the supremum \eqref{pep} is attained for some $\beta>0$ together with $(f,g,A,B,b)$ and others.
In case $\beta \neq 1$, we let $(\bar f,\bar g, \bar A, \bar B, \bar b) := (\beta f, \beta g, \beta  A, \beta B, \beta b)$
and consider
\begin{eqnarray}\label{prob-3}
  \min\bigl\{\bar f(x)+\bar g(y) \mid \bar Ax + \bar B y = \bar b, \;   x\in\R^{n_1}, \;  y\in\R^{n_2}\bigr\},
\end{eqnarray}
which is clearly equivalent to \eqref{prob}. Note that applying the ADMM   to \eqref{prob} with $\beta>0$  is equivalent to applying the ADMM to \eqref{prob-3} with $\beta=1$.
Furthermore, $\Omega^*$ given in \eqref{kkt} can be rewritten as
\begin{equation*}
\Omega^* = \{(x^*, y^*, z^*) \mid \bar A^Tz^* \in \partial \bar f(x^*), \; \bar B^T z^* \in \partial \bar g(y^*), \; \bar Ax^*+ \bar By^*= \bar b\}.
\end{equation*}
Finally, $\beta$ can also be absorbed in both the constraint $R^k=1$ and the objective function $R^{k+1}$, e.g., the constraint  function  $R^k$ can be represented as
\begin{equation*}
R^k =  \|z^k - z^*\|^2 + \gamma \|\bar B(y^k - y^*)\|^2 + (\gamma-1) \|\bar Ax^k+ \bar By^k-\bar b\|^2.
\end{equation*}
This implies that the supremum in \eqref{pep} is also attained for $\beta = 1$ together with  $(\bar f, \bar g, \bar A, \bar B, \bar b)$ and others.
Thus, we may restrict $\beta$ to $1$ without affecting the supremum.

\subsection{Transform $x^*$, $y^*$ and $b$ to $0$}
Let $(x^*,y^*)$ be an optimal solution of \eqref{prob}. Since $Ax^*+By^*=b$, by transformations of variables $\bar x = x - x^*$ and $\bar y = y - y^*$, we can represent \eqref{prob} as
\begin{align}\label{prob-xyb}
\min\bigl\{\bar f(\bar{x})+ \bar g(\bar{y}) \mid A \bar{x} &+B \bar{y} = 0, \; \bar x\in\R^{n_1}, \; \bar y\in\R^{n_2}\bigr\},
\end{align}
where $\bar{f}(\bar x) := f(\bar{x}+x^*) = f(x)$ and $\bar{g}(\bar y) := g(\bar{y}+y^*) = g(y)$.
Similarly, we can apply the same change of variables to the ADMM scheme \eqref{admm2}, the set of KKT points $\Omega^*$, the constraint $R^k=1$ and the objective function $R^{k+1}$. 
Now, the reformulated problem \eqref{prob-xyb} has an optimal solution $(\bar x^*, \bar y^*) = (0,0)$.
Therefore, without affecting the supremum, in the rest of this paper we may always assume $(x^*,y^*) = (0,0)$ and $b=0$.

Recall that we already assumed $\beta=1$.
Then, problem \eqref{prob}, the set of KKT points $\Omega^*$ in \eqref{kkt} and the constraint function $R^k$ in \eqref{Rk} can be simplified to
\begin{align}\label{omega-Rk2}
\left\{
\begin{array}{l}
  \min\bigl\{f(x)+g(y) \mid Ax + By = 0, \;  x\in\R^{n_1}, \; y\in\R^{n_2}\bigr\}, \\[0.2em]
  \Omega^* = \{(0 ,0,z^*) \mid A^Tz^* \in \partial  f(0),\; B^Tz^* \in \partial  g(0)\}, \\[0.2em]
  R^k = \|z^k - z^*\|^2 + \gamma \| B y^k\|^2 + (\gamma-1) \|A  x^k+ B y^k\|^2.
  \end{array}
  \right.
\end{align}
Furthermore, \eqref{admm2}  becomes
\begin{equation}\label{admm4}
\left\{
\begin{array}{l}
  A^T[z^k -  (Ax^{k+1}+By^k)] \in  \partial f(x^{k+1}),\\[0.2em]
  B^T[z^k -  (Ax^{k+1}+By^{k+1})]  \in  \partial g(y^{k+1}),\\[0.2em]
  z^{k+1}=z^k-\gamma (Ax^{k+1}+By^{k+1}).
\end{array}
\right.
\end{equation}

\subsection{Simplified PEP}
In each iteration of the ADMM, $x^k$, $y^k$, and $z^k$ are related.
By letting  $k\leftarrow k-1$ in \eqref{admm4} and eliminating $z^{k-1}$, we obtain
\begin{equation*}
B^T  [z^k+(\gamma-1)(Ax^{k}+By^{k})] \in \partial g(y^{k}).
\end{equation*}
With the above discussions, we can equivalently reformulate the PEP \eqref{pep} as
\begin{equation}
  \label{pep2}
\sup \left\{R^{k+1} \left|
  \begin{array}{c}
  f \text{~and~} g \text{~are closed proper convex}, \\[0.2em]
    A^Tz^*\in\partial f(0), \; B^Tz^* \in \partial g(0), \; R^k = 1, \\[0.2em]
  B^T  [z^k+(\gamma-1)(Ax^{k}+By^{k})] \in \partial g(y^{k}), \\[0.2em]
  A^T[z^k - (A x^{k+1}+B y^k)] \in  \partial f( x^{k+1}),\\[0.2em]
  B^T[z^k - (A x^{k+1}+B y^{k+1})]  \in  \partial g( y^{k+1}),\\[0.2em]
  z^{k+1}=z^k-\gamma(A x^{k+1}+B y^{k+1})
  \end{array}
  \right.\right\}.
\end{equation}
Now, the optimization variables are  $(f,g,A,B)$, $z^*$, $(x^k,y^k,z^k)$, and $(x^{k+1},y^{k+1},z^{k+1})$. The equivalence between \eqref{pep} and \eqref{pep2} is in the sense that both problems share the same optimal function value.

\section{SDP reformulation}\label{sc:sdp}
The PEP \eqref{pep2} remains infinite dimensional and nonconvex due to the constraint ``$f$ and $g$ are closed proper  convex". In this section, we reformulate it as a finite dimensional convex SDP by using cyclic monotonicity property of subdifferential operators and semidefinite lifting.

\subsection{Cyclic monotonicity}
 A set $S\subseteq \R^n \times \R^n$
  is said to be cyclic monotone if  for any subset $\{(x_i, g_i): i = 1,2, \ldots, p\}\subseteq S$ it holds that
  \begin{equation*}\label{cyclic}
  \langle x_1 - x_2, g_1\rangle +   \langle x_2 - x_3, g_2\rangle + \ldots   \langle x_{p-1} - x_{p}, g_{p-1}\rangle +   \langle x_p - x_1, g_p\rangle \geq 0.
  \end{equation*}
Let $h$ be an extended real-valued closed proper convex function on $\R^n$.
It is easy to show from subgradient inequalities that any nonempty subset of $\partial h = \{(x,g)\in \R^{n} \times \R^n: g\in\partial h(x)\}$ is cyclic monotone.
On the other side, if $S = \{(x_i, g_i)\}_{i=1}^q \subseteq \R^{n} \times \R^n$ is cyclic monotone, then there exists an extended real-valued closed proper convex function $h$ on $\R^n$ such that $S\subseteq \partial h$, see \cite[Theorem 24.8]{Roc70} and \cite[Theorem 3.4]{LCNS04jca}.
For smooth strongly convex functions, interpolation conditions without using function values were first developed in \cite[Section 3.4]{Taylor17}.
It thus follows that \eqref{pep2} can be equivalently reformulated as
\begin{equation}
  \label{pep3}
\sup \left\{R^{k+1} \left|
  \begin{array}{c}
R^k = 1, \\[0.2em]
S_1 \text{~and~} S_2 \text{~~are cyclic monotone}, \\[0.2em]
  z^{k+1}=z^k-\gamma(A x^{k+1}+B y^{k+1})
  \end{array}
  \right.\right\},
\end{equation}
where
\begin{eqnarray}\label{def:S}
\left\{
\begin{array}{rl}
  S_1 &:= \Bigl\{\bigl(0, A^Tz^*\bigr), \bigl(x^{k+1},  A^T[z^k - (A x^{k+1}+B y^k)]\bigr)\Bigr\}, \\ [0.6em]
  S_2 &:= \left\{
  \begin{array}{c}
  \bigl(0, B^Tz^*\bigr), \\[0.2em]
  \bigl(y^k, B^T  [z^k+(\gamma-1)(Ax^{k}+By^{k})]\bigr),  \\[0.2em]
  \bigl(y^{k+1},  A^T[z^k - (A x^{k+1}+B y^{k+1})]\bigr)
  \end{array}
  \right\}.
\end{array}
  \right.
\end{eqnarray}
Here, the optimization variables are $(A,B)$, $z^*$, $(x^k,y^k,z^k)$, and $(x^{k+1},y^{k+1},z^{k+1})$.
Again, the equivalence between \eqref{pep2} and \eqref{pep3} is in the sense that both problems share the same optimal function value.
Clearly, \eqref{pep3} is finite dimensional, yet still nonconvex.

Next, we write out explicitly the set of inequalities characterizing the cyclic monotonicity of $S_1$ and $S_2$.
Since $S_1$ has only two points, cyclic monotonicity reduces to a single inequality, that is
\begin{equation}\label{ineq-1}
  \langle Ax^{k+1},\ (z^k-z^*) - (Ax^{k+1}+By^k) \rangle \geq 0.
\end{equation}
On the other hand, $S_2$ has three points and its cyclic monotonicity is equivalent to the following five inequalities
\begin{align}\label{ineq2-6}
\left\{
\begin{array}{rcl}
\langle By^{k},\  (z^k-z^*) +(\gamma-1)(Ax^{k}+By^{k}) \rangle&\geq&0,  \\[0.2em]
\langle By^{k+1},\ (z^k -z^*) - (Ax^{k+1}+By^{k+1})\rangle&\geq&0,  \\[0.2em]
\langle By^k-By^{k+1},\  (Ax^{k+1}+By^{k+1}) + (\gamma-1)(Ax^{k}+By^{k})\rangle&\geq&0,  \\[0.2em]
      \langle By^{k+1},\  - (Ax^{k+1}+By^{k+1}) -  (\gamma-1)(Ax^{k}+By^{k}) \rangle & \\[0.2em]
+\langle By^{k},\   (z^k-z^*) + (\gamma-1)(Ax^{k}+By^{k}) \rangle&\geq&0, \\[0.2em]
       \langle By^k,\    (Ax^{k+1}+By^{k+1}) +  (\gamma-1)(Ax^{k}+By^{k}) \rangle&\\[0.2em]
+\langle By^{k+1},\  (z^k-z^*) - (Ax^{k+1}+By^{k+1}) \rangle&\geq&0.
\end{array}
\right.
\end{align}

\subsection{Grammian representation}
The key observation in the reduction of \eqref{pep3} to an SDP is that all the constraints as well as the objective function are linear in terms of the Grammian matrix $X := P^TP \in \mathbb{S}^5_+$ with
\[
P= [Ax^k,\ By^k,\ Ax^{k+1},\ By^{k+1},\ z^k-z^*] \in \R^{m\times 5},
\]
where $\mathbb{S}^5_+$ denotes the set of  real symmetric positive semidefinite matrices of dimension $5$.
Let $e_i\in \R^5$ be the $i$th unit vector, $i = 1,\ldots,5$. For $u,v\in\R^m$, we define
\[
S(u,v) := (uv^T + vu^T)/2.
\]
Then, it is easy to show that \eqref{ineq-1} is equivalent to
\begin{align*}
 \langle A_1, X \rangle \geq 0 \text{~~with~~}  A_1 := S(e_3, e_5 - e_3 - e_2).
\end{align*}
Similarly,  the inequalities in \eqref{ineq2-6} are respectively equivalent to $\langle A_i, X \rangle \geq0$, $i=2,\ldots,6$, where
\begin{align*}
\left\{
\begin{array}{l}
A_2  := S(e_2,  e_5 +(\gamma-1)(e_1+e_2)), \\[0.2em]
A_3  := S(e_4, e_5 - e_3 - e_4), \\[0.2em]
A_4  := S(e_2-e_4,  (e_3+e_4) + (\gamma-1)(e_1+e_2)), \\[0.2em]
A_5  := S(e_4, -(e_3+e_4) -  (\gamma-1)(e_1+e_2)) + S(e_2,  e_5 + (\gamma-1)(e_1+e_2)), \\[0.2em]
A_6  :=  S(e_2, (e_3+e_4) +  (\gamma-1)(e_1+e_2)) + S(e_4,  e_5 -  e_3 - e_4).
\end{array}
\right.
\end{align*}
Furthermore, $R^k=1$ in \eqref{omega-Rk2} can be represented as $\langle A_7, X\rangle = 1$ with
\begin{align*}
A_7 := S(e_5,e_5)  + \gamma  S(e_2,e_2) + (\gamma-1) S(e_1+e_2,e_1+e_2).
\end{align*}
Similarly, we have
\begin{align*}
R^{k+1} & = \|z^{k+1} - z^*\|^2 + \gamma  \|B y^{k+1}\|^2 + (\gamma-1) \|A x^{k+1}+B y^{k+1}\|^2 \\
& = \|z^{k} - z^* -\gamma(Ax^{k+1}+By^{k+1}) \|^2 + \gamma  \|B y^{k+1}\|^2 + (\gamma-1) \|A x^{k+1}+B y^{k+1}\|^2 \\
& = \langle C, X\rangle,
\end{align*}
where $\langle C, X\rangle := \text{tr}(C^TX)$ denotes the trace inner product and
\begin{align*}
C := S(e_5 - \gamma (e_3 + e_4), e_5 - \gamma (e_3 + e_4)) + \gamma S(e_4, e_4) + (\gamma-1) S(e_3 + e_4,e_3 + e_4).
\end{align*}

\subsection{SDP formulation}
Now, we are ready to represent \eqref{pep3} as the following SDP with rank constraint
\begin{equation}
  \label{pep4}
\max \left\{ \langle C, X\rangle \left|
  \begin{array}{l}
\langle A_i, X\rangle \geq 0, \;  i = 1,\ldots, 6, \\[0.2em]
\langle A_7, X\rangle = 1, \;   X\in \mathbb{S}^5_+, \; \text{rank}(X) \leq m
  \end{array}
  \right.\right\}.
\end{equation}
Again, the equivalence between \eqref{pep3} and \eqref{pep4} are in the sense that they share the same optimal function value.
Since $X\in \mathbb{S}^5_+$, the rank constraint  can be removed safely when $m\geq 5$.
In this case,  \eqref{pep4} is equivalent to the following standard SDP
\begin{equation}
  \label{sdp}
\max \left\{ \langle C, X\rangle  \mid
\langle A_i, X\rangle \geq 0, \;  i = 1,\ldots, 6, \;
\langle A_7, X\rangle = 1, \; X\in \mathbb{S}^5_+
\right\}.
\end{equation}
Now, for $m\geq 5$, we have   reformulated the infinite dimensional nonconvex PEP \eqref{pep} equivalently as the finite dimensional convex SDP \eqref{sdp}.

\begin{figure}[!h]
\centering\begin{tikzpicture}[scale=1.3]
\begin{axis}[xmin=1.5, xmax=2, xlabel={$\gamma$}, ylabel={Optimal function values of SDP}]
\addplot[mark=, blue] table[x=gamma,y=p]{result.txt};
\end{axis}
\end{tikzpicture}
\caption{Optimal function values of \eqref{sdp} for $\gamma \in [1.5,2]$.}
\label{fig}
\end{figure}
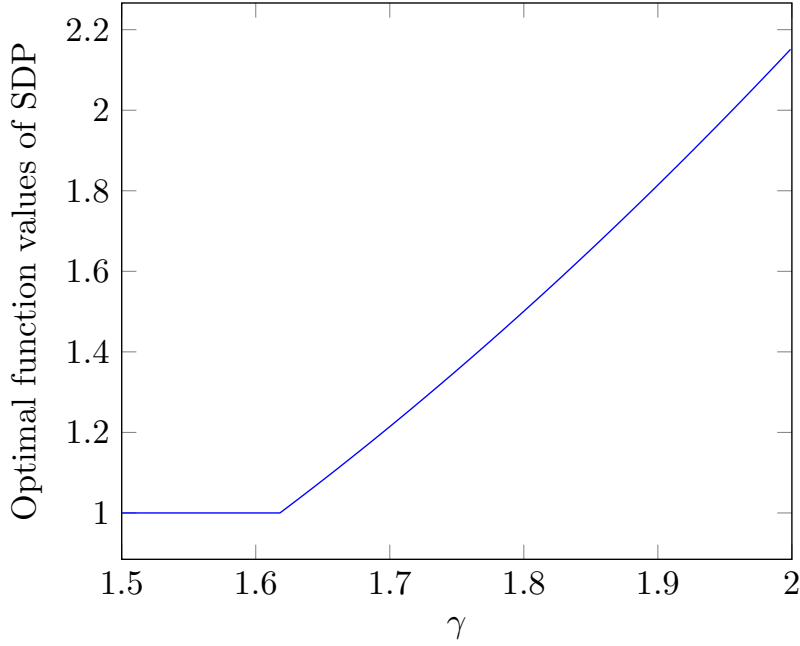

We solved the SDP \eqref{sdp} by SeDuMi \cite{Stu99oms} for $\gamma\in[1.5,2]$. The optimal function values are given in Figure \ref{fig}.
It can be seen from Figure \ref{fig} that there exists a threshold $\phi$ such that the optimal function values are approximately equal to $1$ for $\gamma < \phi$.
This is reasonable because theoretically it is guaranteed that $R^{k+1}\leq R^k$, see, e.g., \cite[Theorem 1]{HY98orl}, and meanwhile linear convergence is generally out of reach in the current setting.
When $\gamma$ gets greater than the threshold $\phi$, the optimal function values become greater than $1$, indicating failure, in the worst case, of the contractiveness of $R^k$ with respect to $k$. Later, we will show that $\phi$ is exactly the golden ratio $(\sqrt{5}+1)/2$.

\section{A rank two feasible solution}\label{sc:rank2}
For $\gamma$ larger than the golden ratio, the numerical solutions of SDP \eqref{sdp} suggests that the constraints are all active and the optimal solution $X$ is of rank 2.
By considering the primal and the dual problems together and assisted by numerical and symbolic computations, we managed to construct a rank two analytic feasible solution to the SDP \eqref{sdp}. It is thus apparent that the constructed solution is also feasible for \eqref{pep4}  in the case $m\geq 2$. Moreover, the objective function value corresponding to the constructed rank two solution is greater than $1$ for $\gamma > (\sqrt{5}+1)/2$.

The constructed rank two feasible solution is given by
\begin{equation}
\begin{cases}
X_{f} = P^TP = \alpha \bar{P}^T\bar{P},\
\alpha=\frac{\gamma^2-2}{-8-\gamma+6 \gamma^2-\gamma^3+(-6+3\gamma+\gamma^2) \sqrt{\gamma^2-1}},\\
\bar{P}(1, 1)=\frac{\bigl(-1-3\gamma+2\gamma^2+(3-2\gamma)\sqrt{\gamma^2-1}\bigr)\sqrt{1+\gamma-\gamma^2+(-1+\gamma)\sqrt{\gamma^2-1}}}{-2-\gamma\bigl(-1-3\gamma+2\gamma^2+(3-2\gamma)\sqrt{\gamma^2-1}\bigr)},\\
\bar{P}(1, 2)=\frac{\sqrt{1+\gamma-\gamma^2+(-1+\gamma)\sqrt{\gamma^2-1}}}{\gamma},\\
\bar{P}(1, 3)=\frac{-\sqrt{\gamma^2-1}}{(1+\gamma)\sqrt{1+\gamma-\gamma^2+(-1+\gamma)\sqrt{\gamma^2-1}}},\\
\bar{P}(1,4)=\frac{\bigl(1-\gamma^2-\gamma\sqrt{\gamma^2-1}\bigr)\sqrt{1+\gamma-\gamma^2+(-1+\gamma)\sqrt{\gamma^2-1}}}{\gamma+\gamma^2},\\
\bar{P}(1,5)=-\frac{2(-1+\gamma)\sqrt{1+\gamma-\gamma^2+(-1+\gamma)\sqrt{\gamma^2-1}}}{\gamma\bigl[2+\gamma\bigl(-1-3\gamma+2\gamma^2+(3-2\gamma)\sqrt{\gamma^2-1}\bigr)\bigr]},\\
\bar{P}(2,1)=\bar{P}(2,2)=\bar{P}(2,4)=0,\\
\bar{P}(2,3)=\frac{1+\gamma+\sqrt{\gamma^2-1}}{1+\gamma},\\
\bar{P}(2,5)=1.
\end{cases}
\end{equation}
The theoretical verification of the feasibility of $X_f$ is tedious, though possible.
Alternatively, we have verified its feasibility via checking \eqref{ineq-1}-\eqref{ineq2-6} and $R^k = 1$ using Mathematica.
Numerically, for the dual step length $\gamma$ greater than the golden ratio,
the objective function value for the constructed rank two feasible solution matches the optimal value of the SDP,
with difference less than $10^{-7}$, implying that the feasible solution constructed is indeed an optimal solution.
We note that $\alpha>0$ for $\gamma>1$ since $R^k = 1 = \alpha \langle A_7, \bar{P}^T\bar{P}\rangle$ and $\langle A_7, \bar{P}^T\bar{P}\rangle$ is a sum-of-squares.

Again, assisted by symbolic computation using Mathematica, we have
\begin{align}
\nonumber
\langle C, X_{f}\rangle
&=  \Biggl[3+\frac{2}{\gamma-1}-2 \gamma^2+\frac{2(-1-\gamma+\gamma^2)\sqrt{\gamma^2-1}}{\gamma-1} \Biggr]^{-1} \\
\label{dedom}
&=  \Biggl[\frac{2(1+\gamma-\gamma^2)}{\bigl(\gamma-1)(\gamma+\sqrt{\gamma^2-1} \bigr)}+1\Biggr]^{-1}.
\end{align}
Since $\text{rank}(X_f)=2$, $X_f$ is also feasible for \eqref{pep4} in the case $m\geq 2$. It thus follows from \eqref{pep4}, \eqref{sdp} and \eqref{dedom} that
\begin{equation*}
\text{val}_{sdp} \geq \text{val}_{pep} \geq \langle C, X_{f}\rangle  > 1 \text{~~for~~} \gamma > \frac{\sqrt{5}+1}{2} \text{~~and~~} m\geq 2.
\end{equation*}
Here val$_{sdp}$ and val$_{pep}$ denote the optimal function values of \eqref{sdp} and \eqref{pep4}, respectively.
This shows that the dual step length $\gamma$  cannot exceed the golden ratio within this analytic framework, at least for problems with $m\geq 2$.

Finally, we note that the equivalence between \eqref{pep} and \eqref{sdp} are bidirectional in the sense that given a feasible solution to one of them, a feasible solution to the other can be constructed accordingly.
On the one hand, it is clear that, given a feasible solution to \eqref{pep}, a matrix $X\in \mathbb{S}^5$ feasible for \eqref{sdp} can be constructed by following the reductions in section \ref{sc:pep} and the transformations in section \ref{sc:sdp}. On the other hand, given $X\in\mathbb{S}^5$ feasible for \eqref{sdp}, one can apply Cholesky factorization to obtain $P$ satisfying $X = P^TP$. Then, the vectors $\{Ax^k, By^k, Ax^{k+1}, By^{k+1}, z^k-z^*\}$ can be extracted form $P$ accordingly. By restricting $A$ and $B$ to be the identity matrices and arbitrarily fixing $z^*$, e.g., at the origin, we obtain $\{x^k, y^k, x^{k+1}, y^{k+1}, z^k, z^*\}$, which, together with $A$ and $B$ being the identity matrices, satisfies the interpolation conditions $S_1 \subseteq \partial f$ and $S_2 \subseteq \partial g$ for some extended real valued closed proper convex function $f$ and $g$, where $S_1$ and $S_2$ are defined in \eqref{def:S}. The functions $f$ and $g$ can then be obtained via standard convex interpolation for cyclically monotone sets, see \cite[Theorem 24.8]{Roc70} and \cite[Theorem 3.4]{LCNS04jca}.

\section{Concluding remarks}\label{sc:remarks}
The performance estimation framework is used to analyzing the ADMM.
Assisted by numerical and symbolic computations, we show that in the general convex setting it is impossible,  at least for $m\geq 2$,
to extend the dual step size $\gamma$ beyond the golden ratio and meanwhile keep the conventional measure $R^k$ decreasing throughout the iterative process.
In our constructed example, we can only guarantee that the conventional measure increases for only one step.
Similar techniques may be used to construct examples in which the conventional measure increases for several consecutive steps,
but the number of cyclic monotone constraints in SDP will blow up factorially.

To the best of our knowledge, it is the first time the performance estimation framework is used to analyze the ADMM-type methods.
Very recently, performance of first order methods using interpolation conditions without function values has been analyzed in \cite{LS20arxiv}, where the difficulty caused by the factorial number of constraints has been alleviated enormously.
However, the analysis in \cite{LS20arxiv} only applies to smooth and strongly convex functions.
%


\def\cprime{$'$}


\newpage
\includegraphics[width=\textwidth]{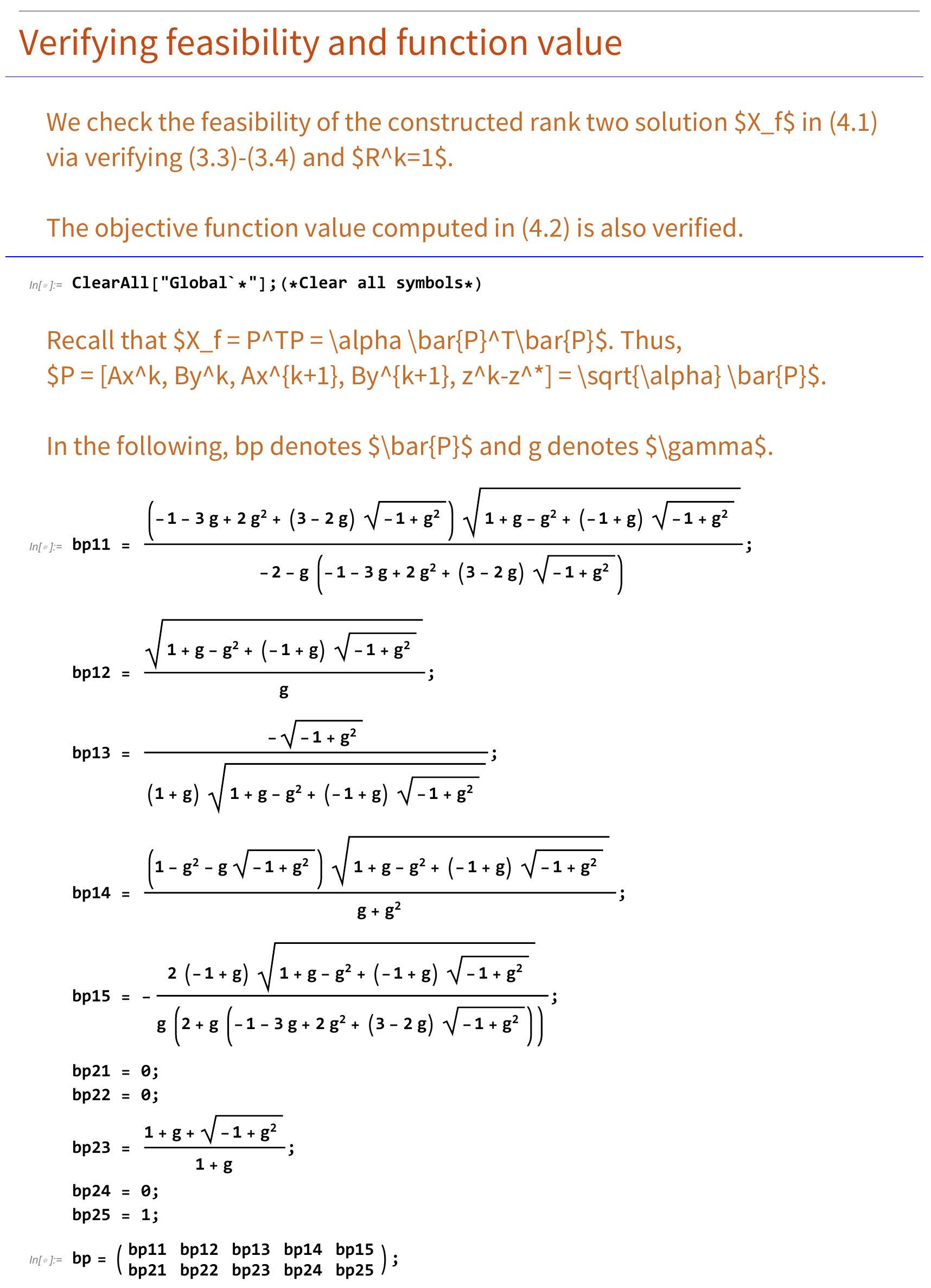}
\includegraphics[width=\textwidth]{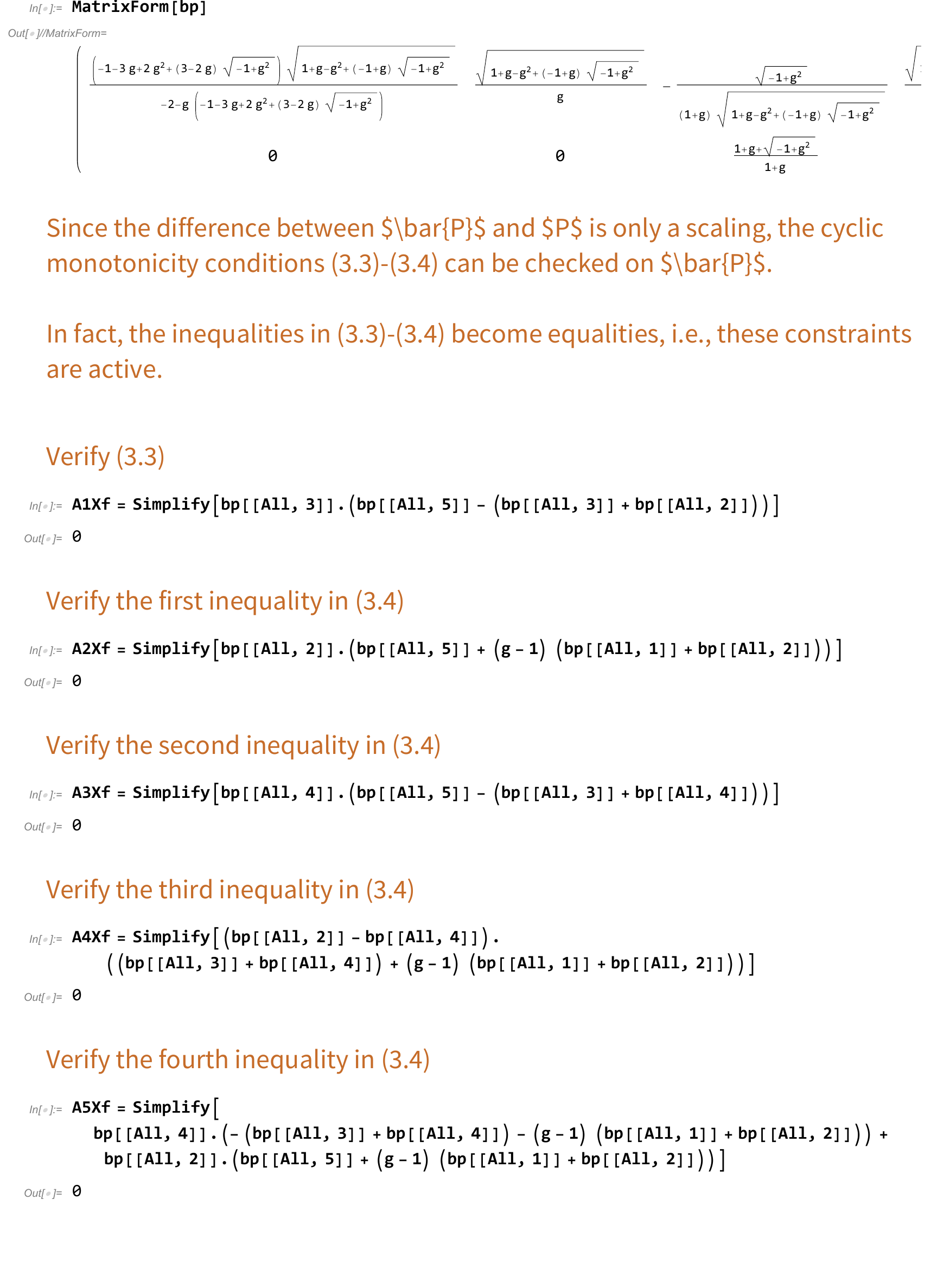}
\includegraphics[width=\textwidth]{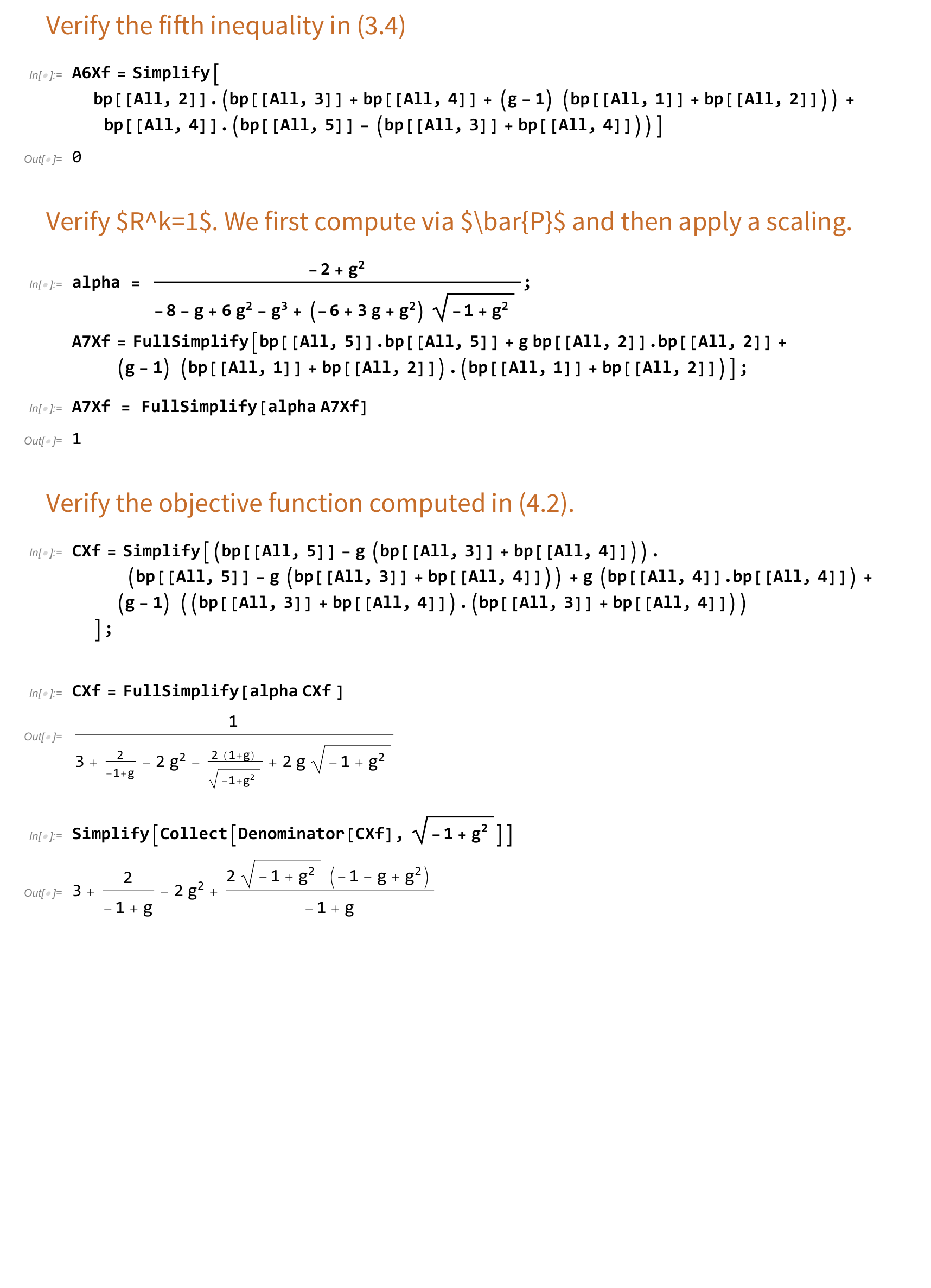}

\end{document}